\documentclass[11pt,onecolumn]{article}
\usepackage{amscd}
\usepackage{times}
\usepackage{amsmath}
\usepackage{amssymb}
\usepackage{xspace}
\usepackage{theorem}
\usepackage{graphicx}
\usepackage{ifpdf}
\usepackage{url,hyperref}
\usepackage{latexsym}
\usepackage{euscript}
\usepackage{xspace}
\usepackage[all]{xy}
\usepackage{color}
\usepackage{makeidx}
\usepackage{tikz}
\long\def\remove#1{}
\newtheorem{theorem}{Theorem}[section] 

\newcommand {\mm}[1] {\ifmmode{#1}\else{\mbox{\(#1\)}}\fi}

\newcommand{\supp}{\mathrm supp}

\newcommand{\rank}                {\mm {\rm rank}}

\newcommand{\cancel}[1]

\begin{document}

\title {Barcodes in level and sub-level persistence and Morse-Novikov theory}

\author{
Dan Burghelea  \thanks{
Department of Mathematics,
The Ohio State University, Columbus, OH 43210,USA.
Email: {\tt burghele@math.ohio-state.edu}}
}
\date{}

\date{}
\maketitle

\begin{abstract}

In this note we recall the relations between the barcodes in level and sub-level persistence and make precise their relation with the Morse-Novikov complex of a Morse real- or angle-valued map. The results in this papers are implicit in \cite{BH} and explicit in \cite {B}, but apparently not well known even to experts.
\end{abstract}

\maketitle

\setcounter{tocdepth}{1}

\vskip .2in

In this note we recall the relations between the barcodes in level and sub-level persistence and make precise their relation with the Morse-Novikov complex of a Morse real- or angle-valued map.  The results in this Note are implicit in \cite{BH} and explicit in \cite {B}. This note  should be viewed as a preliminary section to \cite {B2}.

A short summary of  note is provided by the following two statements :
\begin{itemize} 
\item  From the perspective of classical (ELZ) persistence 
 for a Morse real-valued map  which is Lyapunov for a Morse-Smale vector field, the infinite barcodes determine the Betti numbers of the underlying manifold and the finite barcodes 
give information about the multitude of instantons (visible trajectories) between the rest points. 
\item From the perspective of level persistence  
 for a real-valued or an angle-valued Morse map which is Lyapunov for a Morse - Smale vector field, the closed and the  open barcodes determine the standard  Betti numbers 
 (for real-valued map) resp. the Novikov Betti numbers  (for angle-valued map) while the closed-open barcodes give indication about the multitude of instantons.

More precisely the cardinality of the $r-$closed-open bar codes equals to $N$ indicates that from at least $N$ different rest points of Morse index $r$ originate instantons and to at  least $N $ different rest points  of Morse index $r-1$ arrive instantons.
\end{itemize}

\vskip .1in

\section {Barcodes in level and sub-level persistence}

{\bf The case of real valued maps}
\vskip .1in

{\bf Level persistence}

For a nonnegative integer $r$, a field $\kappa,$  a real-valued {\it proper tame map} $f:X\to \mathbb R$ with $X$ an ANR  ( and based on homology with coefficients in $\kappa$)   the level persistence provides the following four type of barcodes  c.f \cite {BD11}, \cite {BH} or \cite {B} chapters 5 and 6 :
\begin{enumerate}
\item {\it level persistence closed  $r-$barcodes}, $[a,b], a\leq b, a, b\in \mathbb R,$

\item {\it level persistence open $r-$barcodes}, $(\alpha, \beta), \alpha <\beta,  \alpha, \beta\in \mathbb R,$

\item {\it level persistence closed-open  $r-$barcodes}, $[m,n), m<n, m,n \in \mathbb R,$

\item {\it level persistence open-closed  $r-$barcodes}, $(m', n'], m'< n', m', n'\in \mathbb R.$
\end{enumerate}

Recall that $f$ is  {\it proper} if for any compact interval $I\subset \mathbb R,$  $f^{-1}(I)$ is compact, and {\it tame}  if $f^{-1}(t)$ is an  ANR (cf. \cite{B} Chapter 1 for definition) any $t\in \mathbb R.$ 
Recall that all simplicial complexes and in particular all  manifolds, or more general stratified spaces, are ANR's. 

Each such barcode has a multiplicity $\geq 1$ and the collection of barcodes can be recored as  the nonnegative integer-valued maps 
$\delta^f_r : \mathbb C\to \mathbb Z_{\geq 0}$ and $\gamma^f_r : \mathbb C\setminus \Delta \to \mathbb Z_{\geq 0}$

defined as follows:

\begin{equation}
\delta^f_r(z):=\begin{cases} \text{the multiplicity of the r- closed barcode}\  [a,b] \ \text{for} \ z=a+ib\\
\ \text{the multiplicity of the (r-1)- open barcode}\ (\alpha, \beta) \ \text{for} \ z=\beta+i \alpha\\ 
0 \ \text{otherwise},
\end{cases}
\end{equation}

\begin{equation}
\gamma^f_r(z):=\begin{cases} \text{the multiplicity of the r-closed-open  barcode}\  [m,n) \ \text{for} \ z=m+i n\\
\ \text{the multiplicity of the r- open-closed  barcode}\ (m', n'] \ \text{for} \ z=n'+i m'\\ 
0 \ \text{otherwise}
\end{cases}
\end{equation}
cf \cite {BH} or \cite {B}, whose elements are the points in the support of these with multiplicity the value of these maps

If $X$ is compact then the maps $\delta^f_r$ and $\gamma^f_r$ are configurations of points, i.e. maps with finite support.

One convenes  to denote by $\mathcal B_r^c(f),$ $\mathcal B_r^o(f),$ $\mathcal B_r^{c,o}(f),$  $\mathcal B_r^{o,c}(f)$ the multi-sets (sets of elements with multiplicity) of level persistence 
closed $r-$barcodes, open $r-$barcodes, closed-open $r-$barcodes  and $r-$open-closed barcodes.

\vskip .1in
{\bf Sub-level persistence}

For a nonnegative integer $r$, a field $\kappa,$  a real-valued tame proper map $f:X\to \mathbb R$  bounded from below, with $X$ an ANR,  (and based on homology with coefficients in $\kappa$)  the sub-level persistence provides the following two type of barcodes  c.f \cite {ELZ}
\begin{enumerate}
\item {\it sub-level persistence infinite $r-$barcodes} $(a,\infty),  a\in \mathbb R,$

\item {\it sub-level persistence finite  $r-$barcodes} $(a,b), a < b, \  a,b \in \mathbb R.$
\end{enumerate}

One convenes  to denote by $S\mathcal B_r^\infty(f)$ and $S\mathcal B_r^{finite}(f)$  the multi-sets  of the sub-level  persistence 
infinite $r-$barcodes and finite $r-$barcodes respectively.

Level persistence barcodes refine the sub-level persistence barcodes as follows:

\begin{itemize}
\item A level persistence closed $r-$barcode $[a,b]$ corresponds to a sub-level persistence  infinite $r-$barcode $(a,\infty)$ with the same multiplicity.

\item A level persistence open $r-$barcode $(\alpha, \beta)$ corresponds to a sub-level persistence  infinite $(r+1)-$barcode $(\beta,\infty)$ with the same multiplicity.

\item A level persistence closed-open $r-$barcode $[m, n)$ corresponds to a sub-level persistence  finite $r-$barcode $(m,n)$ with the same multiplicity.
\end{itemize}
The level persistence open-closed  $r-$barcode $(m', n']$ is not visible in sub-level persistence of $f$ (but it does appear in the sub-level persistence of $-f.$

Note that the ends of any barcode for $f$ are among the critical values of $f$, the values $t$ for which the homology of the level of $f$ at $t$ differs from  the homology of levels at values in an arbitrary neighborhood of $t.$
\vskip .1in

{\bf The case of angle-valued maps}
\vskip .1in
If $f:X\to \mathbb S^1$ is a tame angle-valued map  let  $\tilde f:\tilde X\to \mathbb R$ be a lift of $f$, i.e. given by a  pullback diagram  

$$\xymatrix { \mathbb R\ar[r]^{\pi}&\mathbb S^1\\
\tilde X\ar[r]\ar[u]^{\tilde f}&X\ar[u]^f}$$ with $\pi$ the canonical infinite cyclic cover.
Recall that {\it tame} means $f^{-1}(\theta)$ is an ANR for any $\theta\in \mathbb S^1.$
One defines the multi-sets $\mathcal B^c_r(f),  \mathcal B^o_r(f), \mathcal B^{c,o}_r(f), \mathcal B^{o,c}_r(f)$ as the quotient sets $$\mathcal B^{\cdots}_r( f):= \mathcal B^{\cdots}_r(\tilde f)/ {2\pi \mathbb Z}$$ with respect to the obvious free action of $\mathbb Z$ given by $\mu (k, I)= I+2\pi k.$ The multiplicity of a an element in the quotient set is  the same as the multiplicity of any of its representatives.\vskip .1in

{\bf Barcodes for exact or an integral (topological) closed one form}

If one likes to consider tame real-valued maps $f:X\to \mathbb R$ up to an additive constant, i.e. as an exact topological closed one form $\omega= df,$ 
\footnote {the concept of topological closed one form will be defined in \cite {B2} and is not used in  this note}  or tame  angle-valued maps $f: X\to \mathbb S^1$ up to composition by a rotation of $\mathbb S^1,$  i.e. as a topological integral close one form $\omega= f^\ast (dt)$ with $dt$ the canonical differential one form  on $\mathbb S^1$ (the length) then 
the barcodes of $\omega$ are actually real numbers given by the length  of a barcode $I\in \mathcal B_r^{\cdots}(f)$ or $I\in \mathcal B_r^{\cdots}(\tilde f),$ 
precisely $p(I)= b-a,$  $p:\mathbb C\equiv \mathbb R^2 \to \mathbb R,$  $p(z= a+ib)=b-a,$ with the multiplicity of the barcode $x\in \mathbb R$ given the sum of the multiplicity of each barcode $I\in p^{-1}(x).$

This  explains why the barcodes for an (arbitrary) closed one form which will be introduced in a future paper \cite {B2} are real numbers.

\section{Chain complexes of finite dimensional vector spaces over a field $\kappa$ }
\vskip .1in
{\bf Some elementary linear algebra}
\vskip .1in
It is a straightforward consequence of elementary linear algebra that  a chain complex, cf \cite{B} chapter 8, 

$(C_\ast, \partial_\ast):= \xymatrix{\cdots \ar[r]&C_n \ar[r]^{\partial _n} &C_{n-1}\ar[r]^{\partial_ {n-1}}&\cdots\ar[r]^{\partial _2}&C_0\ar[r]^{\partial _0} &0}$ 

is determined up to a non canonical isomorphism by two of the  three systems of integers:
\begin{enumerate}
\item $c_r:=\dim C_r$
\item $\beta_r:= \dim H_r(C_\ast, \partial _\ast)$
\item $\rho_r:= \rank \ \partial _r$
\end{enumerate}
related  by 
$$c_r= \beta_r +\rho_r +\rho_{r-1}$$
 
 Actually  the chain complex is isomorphic to the chain complex  $(^hC_\ast,^h\partial _r)$ with 
 $^hC_n \simeq  \kappa^{\beta_r}\oplus  \kappa^{\rho_r}\oplus  \kappa^{\rho_{r-1}}$
 
 $^h\partial _r= \begin{pmatrix} 0& 0& 0 \\0& 0& id \\0&0&0\end{pmatrix}.$  This chain complex $(^hC_\ast, ^h\partial_\ast)$is referred to as the {\it  Hodge form} of $(C_\ast, \partial_\ast)$ and is isomorphic although non canonically to $(C_\ast,\partial_\ast).$
\vskip .1in

Suppose that the chain complex $(C'_\ast, \partial '_\ast)$ is a sub-complex of  $(C_\ast, \partial _\ast)$  i.e. $C'_\ast \subseteq C_\ast$ and $\partial '_\ast = \partial _r\mid _{C'_r}$ and therefore $c'_r\leq c_r$ and $\rho'_r\leq \rho_r.$ 

It is straightforward to check that any isomorphism from  $i'_\ast: (C'_\ast, \partial '_\ast)\to (^hC'_\ast, ^h\partial' _\ast)$ can be extended to an isomorphism $i_\ast: (C_\ast, \partial _\ast)\to (^hC_\ast, ^h\partial _\ast).$
\vskip .2in

{\bf Morse complex}
\vskip .1in
Let $M^n$ be a closed manifold, $X$ a smooth vector field which is Morse-Smale and $f:M^n\to \mathbb R$ a Morse function which is Lyapunov for $X$. 
Each rest point $x\in \mathcal X:=\{y\in  M^n\mid X(y)=0\},$ equivalently critical point of $f,$ has a Morse index $r\in \{0, 1,\cdots n\}.$  Let $\kappa$ be a field.

To these data one associate the $\kappa-$vector spaces $C_r:=\kappa[\mathcal X_r]$ generated by the set $\mathcal X_r$ of the rest points of Morse index $r$   
and the linear maps $\partial _r: C_r\to C_{r-1}$ given by the matrix $\partial _r$ with entries $\partial _r(x,y)$ which  counts the algebraic cardinality of the set of 
trajectories of $X$ from the rest point $x\in \mathcal X_r$ to $y\in \mathcal X_{r-1}$ (instantons) viewed as an element in $\kappa.$   

For $t\in \mathbb R$ one denotes by $\mathcal X_r(t):=\{x\in \mathcal X_r, f(x)\leq t\}$ and one observes that the subspace $C_r(t):=\kappa [\mathcal X_r(t)]$ is invariant to $\partial_r ,$ hence $\partial _r(C_r(t))\subseteq C_{r-1}(t).$  
Let $\partial _r(t)$ be the restriction of $\partial _r$ to $C_r(t).$ This equips $(C_\ast, \partial _\ast)$ with the filtration indexed by $t\in \mathbb R,$

$$0
\subseteq  (C_\ast(t'), \partial _r(t'))\subseteq  (C_\ast(t), \partial _r(t))\cdots \subseteq  (C_\ast, \partial _r)$$ for $t'<t.$

One of the main results of Morse  theory is the following theorem.
\begin{theorem}(Morse theory)\

\begin{enumerate}
\item $\partial _{r-1}\cdot \partial _r=0,$ hence $(C_\ast, \partial _\ast)$ defines a chain complex of finite dimensional vector spaces refereed below as the Morse complex of the pair $(M,X)$ and $(C_\ast(t),\partial _\ast (t))$ provides a filtration of the chain complex $(C_\ast, \partial _\ast).$
\item  $H_r(C_\ast,\partial _\ast)\simeq H_r(M;\kappa)$ and $H_r(C_\ast(t),\partial_\ast(t))= H_r(f^{-1}((-\infty, t]), \kappa),$
\item Up to a non canonical isomorphism the Morse complex $(C_\ast,\partial _\ast)$ and its sub-complexes $(C_\ast(t),\partial _\ast(t))$ depend only on the Morse function $f.$ 
\end{enumerate}
\end{theorem}

One can extend the above result to the case $M$ is a compact manifold with boundary $\partial M$ with $X$ transversal to $\partial M$and pointing out towards the interior and $f$ is constant on $\partial M.$ The statement remains the same.

A similar complex is obtained in case $f: X\to \mathbb S^1$ is an angle-valued map  which is Lyapunov for $X$. In this case the set of  trajectories from $x\in \mathcal X_r$ to $y\in \mathcal X_{r-1}$ can have countable cardinality  but their "counting" can still be done as shown by Novikov cf \cite {F} or \cite {P} with the result an element  in the field $\kappa[t^{-1},t]]$ of Laurent power series  with coefficients in $\kappa.$ 

Finally the complex referred below as the Novikov complex has the components $C_k$ the $\kappa[t^{-1},t]]-$vector space generated by $\mathcal X_k$ and the linear maps $\partial _k$ given by matrices with coefficients $\partial _k(x,y), x\in \mathcal X_r, y\in \mathcal X_{r-1},$ elements in $\kappa[t^{-1},t]].$
 
There is no Morse type  filtration  in this case \footnote {but there are other interesting filtrations}.

The above Morse complex and its Morse filtration  as well as the Novikov complex  can be recovered up to a non canonical isomorphism from the barcodes of the level persistence.  More precisely one has.

\vskip .1in

\begin{theorem}\ 

If $f:M^n\to \mathbb R$ is a Morse  real-valued function then $f$ is tame  and  
\begin{enumerate}
\item
$$\beta_r= \sum_{z\in \supp \delta^f_r} \delta^f_r(z)= \sharp (\mathcal B_r^c \sqcup \mathcal B^o_{r-1})= \sharp S\mathcal B_r^{\infty}$$

$$\rho_r= \sum_{z\in \supp \delta^f_r, \Re z <\Im z}\gamma^f_r(z)= \sharp S\mathcal B_r^{finite}= \sharp \mathcal B^{c,o}_r$$

and then 
$$c_r= \beta_r+ \rho_r + \rho_{r-1}.$$
\item
$$\beta_r(t)= \sharp \begin{cases} \{I\in \mathcal B_r^c \mid  I\cap (-\infty,t]\ne \emptyset\} \sqcup \\
\{I\in \mathcal B_r^{c,o} \mid  I\cap (-\infty,t]\ne \emptyset, I\cap [t,\infty)\ne \emptyset\} \sqcup
\\
\{I\in \mathcal B_{r-1}^o \mid  I\subset (-\infty,t]\ne \emptyset\}
\end{cases}$$

$$\rho_r(t)= 
\sharp \{I\in \mathcal B^{c,o}_r\mid I\subset (-\infty, t]\}  $$
where $\sharp$ denotes the cardinality of the multi-set in parentheses.
Then  one has
$$c_r(t)= \beta_r(t)+ \rho_r (t)+ \rho_{r-1}(t).$$
\end{enumerate}
\end{theorem}

In view of the observation that the end of each barcode is a critical value  the knowledge of the $(C_\ast(t), \partial _\ast(t)),$ for all critical values $t,$ determines the entire collection of bar codes $\mathcal B^c_r,$  $\mathcal B^o_r,$ and $\mathcal B^{c,o}_r.$

For  $I$ a barcode with ends $a,b$ $a\leq b$ denote by $l(I)=a$ and $r(I)= b.$ As shown in \cite{B} Chapter 4 one has:
 
\begin{enumerate}
\item
$\mathcal H_r(t):= \begin{cases} \{I\in \mathcal B^c_r \mid l(I)\leq t\}\sqcup\\ \{I\in \mathcal B^o_{r-1} \mid r(I)\leq t\}\sqcup\\  \{I\in \mathcal B^{c,o}_{r} \mid    l(I) \leq t, r(I)> t\}\end{cases}$ with $H_r(C_\ast(t), \partial _\ast(t))\simeq  \kappa[\mathcal H_r(t)],$
\item $\mathcal C^+_r(t):= \{ I\in \mathcal B^{c.o}_r \mid    l(I) \leq t, r(I)\leq t\}$ \ with $C_r^+(t)\simeq \kappa[\mathcal C_r^+(t)],$

\item $\mathcal C^-_r(t):= \{ I\in \mathcal B^{c.o}_{r-1} \mid     r(I)< t\}$ \ with $C_r^-(t)\simeq \kappa[\mathcal C_r^-(t),$
and 
\item $\mathcal C_r(t):= \mathcal H_r(t) \sqcup \mathcal C^+_r(t)\sqcup \mathcal C^-_r(t)\subseteq \mathcal B^c_r\sqcup \mathcal B^o_{r-1}\sqcup \mathcal B^{c,o}_r \sqcup \mathcal B^{c,o}_{r-1}$ \ with $C_r\simeq \kappa[\mathcal C_r].$
\end{enumerate}

For $t\leq t'$  let $i_r(t,t'): \kappa[\mathcal C_r(t)] \to \kappa[ \mathcal C_r(t')]$ be the inclusion induced linear map 
and denote by $\partial _r$ the linear extension of 
 $^h\partial _r(I):= \begin{cases} 0 \ \text{if} \ I\in \mathcal H_r\sqcup \mathcal C_r^+\\ I\  \text{if}\  I\in  \mathcal C_{r-1}^-\end{cases}.$

\begin{theorem}\

Suppose $f: M^n\to \mathbb R$ is a Morse function defined on a closed manifold.
The linear maps $i_r(t,t')$  and $i_{r-1}(t,t')$ intertwine  $^h\partial _r(t)$ and $^h\partial _r(t')$ and then for $t< t',$ the complexes $(^hC_\ast (t),^h\partial_\ast(t))$ provide a filtration of 
$(^hC_\ast ,^h\partial C_\ast)$ which makes the filtered chain complex $(^hC_\ast ,^h\partial _\ast)$ non canonical isomorphic to the Morse complex equipped with a Morse filtration.
\end{theorem}

\end{document}